\documentclass{amsart}
\usepackage{amssymb,amsmath,eucal}

\def\R {{\mathbb R }}

\def\K{{\Bbb K}}
\def\H{{\mathbb H}}

\def\SL{{\rm SL}}
\def\SU{{\rm SU}}

\def\ov{\overline}
\def\phi{\varphi}
\def\epsilon{\varepsilon}
\def\kappa{\varkappa}

\newcounter{sec}

\newcounter{punct}[sec]

\def\punct{\refstepcounter{punct}{\arabic{sec}.\arabic{punct}.  }}

\newtheorem{theorem}{Theorem}[sec]
\newtheorem{proposition}[theorem]{Proposition}

\newtheorem{lemma}[theorem]{Lemma}

\def\Aut{\mathrm{Aut}}

\def\G{\mathbb{G}}
\def\H{\mathbb{H}}

\def\COUNTERS{\addtocounter{sec}{1}
              \setcounter{punct}{0}
          \setcounter{equation}{0}
          \setcounter{theorem}{0}
          \setcounter{problem}{0}
          }
          
          \def\sm{\smallskip}

\begin{document}

\begin{center}
{\bf\Large
Several remarks on  groups of automorphisms of 
free groups}

\bigskip

{\sc \large Yu.A.Neretin}%
\footnote{Supported by grants FWF, P22122, P25142.}
\end{center}

\bigskip

{\small {\sc Abstract.}
Let $\G$ be the group of automorphisms of a free group $F_\infty$
of infinite order. Let $\H$ be the stabilizer of
first $m$ generators of $F_\infty$.
We show that the double cosets 
$\Gamma_m=\H\setminus \G/\H$
admit a natural semigroup structure. For any compact group 
$K$ the semigroup $\Gamma_m$ acts in the space $L^2$ on the 
product of $m$ copies of $K$.}

\bigskip

\section{Statements}

\COUNTERS

{\bf\punct Notation.}
Let $F_n$ be the free group with $n$ generators $x_1$, \dots, $x_n$. Let
$$
F_\infty=\lim_{n\to\infty} F_n
$$
be the free group with countable number of generators. Denote by $\Aut(F_n)$ the group
of automorphisms of $F_n$. Such automorphisms are determined by images
of the generators,
\begin{equation}
x_1\mapsto \gamma_1(x),\quad \dots, \quad x_n\mapsto \gamma_n(x)
,
\label{eq:lll}
\end{equation}
where $\gamma_1(x)$, \dots, $\gamma_n(x)$ is a  collection of elements of $F_n$
(this collection is not arbitrary, for an introduction to the theory of groups
$\Aut_n$, see \cite{LS}, for relatively recent developments,
see survey \cite{Vog}). We use a symbolic notation for (\ref{eq:lll}),
$$
x\mapsto\gamma(x).
$$

 Denote by
$$
\G=\Aut(F_\infty)=\lim_{n\to\infty} \Aut(F_n)
$$
the inductive limit of the groups $\Aut(F_n)$. The group $\Aut(F_\infty)$ acts by
automorphisms on $F_\infty$, but it is not the whole group of automorphisms of
$F_\infty$.
Denote by $\H=\H_m\subset \G$ the stabiliser of the generators $x_1$,\dots $x_m$. 
It is convenient to rename generators as
$x_1$, \dots, $x_m$, $y_1$, $y_2$, \dots
and to write elements of $\H$ as
$$
\begin{cases}
x\mapsto x
\\
y\mapsto \sigma(x,y)
\end{cases},
$$
or, detailed
$$
x_1\mapsto x_1,\quad\dots,\quad x_m\mapsto x_m,\quad y_1\mapsto\sigma_1(x,y),\quad y_2\mapsto\sigma_2(x,y),
\quad \dots
$$

Denote by $S_\infty$ the group of finitely supported permutations of $\{1,2,\dots,\}$.
The  group  $S_\infty$ acts on $F_\infty$ by permutations of generators, this determines an embedding
$S_\infty\to \G$. We also regard $S_\infty$ as a group of infinite 0-1 matrices.
Finally, we define the group
$S_\infty[m]\subset S_\infty$
consisting of substitutions preserving elements $1$, $2$, \dots, $m$.

\sm

{\bf\punct Forcing apart and the semigroup of double cosets.}
Denote by $\Gamma_m$ the double cosets
$$
\Gamma_m:= \H\setminus \G/\H.
$$
Consider the following sequence $\theta_j\in S_\infty[m]$:
$$
\theta_j=
\theta_j[m]:=
\begin{pmatrix}1_m& 0&0&0\\
0&0&1_j&0\\
0&1_j&0&0\\
0&0&0&1_\infty  \end{pmatrix}
.
$$

Fix $g$, $h\in \G$. Consider the following sequence of double
cosets
$$
\H\cdot g \theta_j h\cdot \H
.$$
Evidently, this sequence is eventually constant. We denote
by $g\circ h\in \Gamma_m$ its value for sufficiently large $j$.

\begin{proposition}
\label{pr:Gamma}
{\rm  a)} The double coset containing $g\circ h$ depends only on double cosets containing
$g$ and $h$.

\sm

{\rm b)} The operation $g\circ h$ on $\Gamma_m$ is associative.
\end{proposition}

Let us rename the generators of $F_\infty$, denote them by 
$x_1$, \dots, $x_m$, $y_1$, \dots, $y_N$, $z_1$,\dots, $z_N$, $u_1$, $u_2$, \dots,
where $N$ is sufficiently large (such that both automorphisms $g$, $h$ fix generators $z$, $u$).
Let us write $g$, $\theta_N$, $h$ as
\begin{equation}
\label{eq:product}
g:\,\,
\begin{cases}
x\mapsto \alpha(x,y)\\
y\mapsto \beta(x,y)\\
z\mapsto z\\
u\mapsto u
\end{cases}
\quad
\theta_N:
\begin{cases}
x\mapsto x\\
y\mapsto z\\
z\mapsto y\\
u\mapsto u
\end{cases}
\quad
h:\,\,
\begin{cases}
x\mapsto \gamma(x,y)\\
y\mapsto \delta(x,y)\\
z\mapsto z\\
u\mapsto u
\end{cases}
\end{equation}
Then the product is 
\begin{equation}
g\circ h
:\begin{cases}
x\mapsto \gamma(\alpha(x,y),z)
\\
y\mapsto \delta(\alpha(x,y),z)
\\
z\mapsto \beta(x,y)
\\
u\mapsto u
\end{cases}
\label{eq:product2}
\end{equation}

{\sc Remark.}
The group of invertible elements of $\Gamma_m$ is $\Aut(F_m)$.
\hfill $\square$

\sm

 {\bf\punct Actions in $L^2$.} Let $K$ be a compact group, $U\subset K$ a closed subgroup%
 \footnote{The main interesting case is $K=U=\SU(2)$.}.
 Equip $K$ with a probabilistic Haar measure.
 Consider the countable product 
 $$
 K^\infty:=K\times K\times \dots
 $$
The group $U$ acts on $K^\infty$ by conjugations
 \begin{equation}
(k_1,k_2,\dots)\mapsto (uk_1u^{-1},uk_2u^{-1},\dots)
\label{eq:uuu}
\end{equation}
 Consider the space of {\it conjugacy classes} $K^\infty//U$.

Let $k=(k_1,k_2,\dots)\in K^\infty$. For any element  
(\ref{eq:lll}) of  $\G$ we define the map
$K^\infty\to K^\infty$ given by
\begin{equation}
g:(k_1, k_2, \dots)\mapsto
\bigl(l_1(k_1, k_2, \dots), l_2(k_1, k_2, \dots), \dots\bigr)
\label{eq:kkk}
\end{equation}
(we substitute $k_1$, $k_2$, \dots to the corresponding words).
Thus we get an action of the group $\G$ on the space
$K^\infty$. Such maps preserve the Haar measure on $K^\infty$
(this is clear for generators of $\G$, on presentation
of this group,  see, e.g., \cite{LS}, Section 1.4). 

The transformations (\ref{eq:kkk}) commute with the action
(\ref{eq:uuu}) of $U$. Therefore we get a measure preserving action of 
$\G$ on $K^\infty//U$ and 
the unitary representation 
$$
T(g)f(k)=f(g(k))
$$
of $\G$
in $L^2(K^\infty//U)$.

Denote by $H\subset L^2(K^\infty)$ the space of functions depending 
only on $k_1$,\dots, $k_m$,
$$H\simeq L^2(K^m).$$
 Denote by $P$ the operator of orthogonal projection
to $H$. For $g\in\Aut_\infty$ we define  the operator
$$
\ov T(g): H\to H
$$
given by
$$
\ov T(g)=P\,T(g).
$$ 
Evidently, for any $h_1$, $h_2\in\H$ we have
$$
\ov T(k_1gk_2)=\ov T(g)
.
$$
Hence $g\mapsto \ov T(g)$ is a well-defined operator-valued function 
on the semigroup $\Gamma_m$

\begin{theorem}
\label{th}
$\ov T$ is a representation of the semigroup $\Gamma_m$ in $L^2(K^m)$.
\end{theorem}

{\sc Remark.} This also determines an action of $\Gamma_m$ on 
the measure space $K^m//U$ by polymorphisms (spreading maps, see
\cite{Ner-book}, Section VIII.4.
\hfill $\square$

\sm

{\bf\punct Some comments.}
The phenomenas discussed above (the existence of semigroup structure of double cosets
and the action of the semigroup
in the subspace of fixed vectors) are usual for infinite-dimensional groups. First special cases were
discovered by R.S.Ismagilov in 60s (see \cite{Ism1}, \cite{Ism2}). The phenomenas exist for classical groups
over $\R$ and over $p$-adic fields, for symmetric groups, for groups
of automorphisms of measure spaces. This was widely explored by G.I.Olshanski 
in representation theory of infinite-dimensional classical groups (see \cite{Olsh-GB}, \cite{Olsh-topics}, see also
\cite{Ner-book}). On recent progress, see, e.g., \cite{Ner-symm}, \cite{Ner-char}, \cite{Ner-faa}.
The present note shows that a behavior of $\Aut_\infty$ is (at least partially)
 similar to the behavior of infinite-dimensional groups. The proof of Theorem \ref{th}
 given below (\ref{eq:long}) coincides with a proof of \cite{Ner-book}, Theorem VIII.5.1.

The spaces of conjugacy classes $(K\times\dots\times K)//K$ and actions of discrete groups
on these spaces are widely discussed in theory
  of Teichm\"uller spaces and its neighborhood (on actions in $L^2$, see \cite{Gol}, \cite{Gol1}, \cite{PX}).
  
\sm

{\bf\punct Extensions of the construction.}
Consider the conjugacy classes 
$$
\Delta_m:=\G//\H.
$$
For $g$, $h\in\H$  we consider the conjugacy class
containing
$$
g\,  (\theta_j h \theta_j^{-1})
.
$$
This sequence is eventually constant, we set $g\ast h$ being its value for large $j$.

\begin{proposition}
\label{pr:conjugacy}
 The $\ast$-multiplication is a well-defined associative operation
 on $\Delta_m$.
\end{proposition}

Next, consider a product $\G^k=\G\times \dots \times\G$
of $k$ copies of $\G$. Consider the diagonal subgroup $\G=\mathrm{diag}(G)$
and the subgroup $H\subset$ in the diagonal.

For $g$, $h\in \G$ consider the following  sequence of double cosets
$$
\H\cdot g \theta_j h\cdot \H
,$$
where $\theta_j$ is regarded as an element of $\H$.
Again, this sequence is eventually constant, we denote by $g\circ h$ its value for
sufficiently large $j$

\begin{proposition}
\label{pr:mnogo}
The $\circ$-multiplication is a well-defined associative operation on $\H\setminus \G^k/\H$.
\end{proposition}

{\bf Acknowledgments.} I am grateful to P.Michor for discussion of this topic. 

\section{Proofs}

\COUNTERS

{\bf\punct Proof of Proposition \ref{pr:Gamma}.}
First, we show that  the product does not depend on the choice of $N$. Indeed, let us denote generators as
$$
x_1,\dots,x_m, y_1,\dots,y_N,y'_1,\dots,y'_p, z_1,\dots,z_N,z'_1,\dots,z'_p, u_1, u_2,\dots
$$
We get
$$
g\theta_{N+p} h:
\begin{cases}
x\mapsto \gamma(\alpha(x,y),z)
\\
y\mapsto \delta(\alpha(x,y),z)
\\
y'\mapsto z'
\\
z\mapsto \beta(x,y)
\\
z'\mapsto y'
\\
u\mapsto u
\end{cases}
$$
We multiply this automorphism by a substitution
$$
y'\mapsto z',\qquad z'\mapsto y' 
$$
and renumerate generators in the order
$$
x_1,\dots,x_m, y_1,\dots,y_N,z_1,\dots,z_N,y'_1,\dots,y'_p, z'_1,\dots,z'_p, u_1, u_2,\dots
$$       
Then we obtain (\ref{eq:product2}).       
Such renumeration is equivalent to a conjugation of $g\theta_{N+p} h$ by a certain 
element of $S_\infty[m]$.

Next, consider  elements $r$, $q\in\H$
given by 
$$r:
\begin{cases}
x\mapsto x\\
y\mapsto \sigma(x,y)\\
z\mapsto z\\
u\mapsto u
\end{cases}
\qquad
q:
\begin{cases}
x\mapsto x\\
y\mapsto \tau(x,y)\\
z\mapsto z\\
u\mapsto u
\end{cases}
$$
Then
$$
g\theta_N r h:
\begin{cases}
x\mapsto \gamma(\alpha(x,y),\sigma(\alpha(x,y), z))
\\
y\mapsto \delta(\alpha(x,y),\sigma(\alpha(x,y),z))
\\
z\mapsto \beta(x,y)
\\
u\mapsto u
\end{cases}
$$
Therefore,
$$
g\theta_N r h=
r^\square\cdot g\theta_N  h
,
$$
where $r^\square$ is an endomorphism of $F_\infty$ given by
$$
r^\square:
\begin{cases}
 x\mapsto x\\
  y\mapsto y\\
  z\mapsto \sigma(\alpha(x,y),z)\\
   u\mapsto u\\
\end{cases}
$$
To show invertibility of $r^\square$,
we write $r^{-1}$ as
$$r^{-1}:
\begin{cases}
x\mapsto x\\
y\mapsto s(x,y)\\
z\mapsto z\\
u\mapsto u
\end{cases}
.$$
Then
$$
s(\sigma(x,y),y)=y\qquad \sigma(s(x,y),y)=y
$$
In these equations we can replace $x_1$,\dots, $x_m$ by an arbitrary collection
of words in $F_\infty$ without entries of $y$, Therefore,
the endomorphism
$$
\begin{cases}
 x\mapsto x\\
  y\mapsto y\\
  z\mapsto s(\alpha(x,y),z)\\
   u\mapsto u\\
\end{cases}
$$
is inverse to $r^\square$.

Next, consider $g q \theta_N h$. Passing to the inverse element
$$
(g q \theta_N h)^{-1}=h^{-1} \theta_N q^{-1} g^{-1},
$$
we come to the the case  discussed just now,
$$
h^{-1} \theta_N q^{-1} g^{-1}= q^\triangledown\cdot h^{-1} \theta_N  g^{-1}
$$
for some $q^\triangledown\in \K$.
Therefore,
$$
g q \theta_N h=g \theta_N h \cdot (q^\triangledown)^{-1}.
$$

 This proves the statement a).

\sm

To prove associativity,
take 3 elements of $\G$,
$$
g:\,\,
\begin{cases}
x\mapsto \alpha(x,y)\\
y\mapsto \beta(x,y)\\
z\mapsto z
\end{cases}
\quad
h:\,\,
\begin{cases}
x\mapsto \gamma(x,y)\\
y\mapsto \delta(x,y)\\
z\mapsto z
\end{cases}
f:
\begin{cases}
x\mapsto \phi(x,y)\\
y\mapsto \psi(x,y)\\
z\mapsto z
\end{cases}
$$
 To evaluate 
\begin{equation} 
\label{eq:ass}
 (g\circ h)\circ f,\qquad g\circ(h\circ f)
 \end{equation}
 we can change these elements by  conjugate elements
 $$
g:\,\,
\begin{cases}
x\mapsto \alpha(x,y)\\
y\mapsto y \\
z\mapsto z\\
u\mapsto \beta(x,u)\\
v\mapsto v
\end{cases}
\quad
h:\,\,
\begin{cases}
x\mapsto \gamma(x,y)\\
y\mapsto y\\ 
z\mapsto \delta(x,z)\\
u\mapsto u\\
v\mapsto v
\end{cases}
f:
\begin{cases}
x\mapsto \phi(x,y)\\
y\mapsto \psi(x,y) \\
z\mapsto z\\
u\mapsto u\\
v\mapsto v
\end{cases}
$$
 Then for calculation of $\circ$-products in ({\ref{eq:ass}) we can  set
 $\theta_N=\theta_0$, i.e., we can evaluate the usual product.
 Now associativity is obvious. The final formula is
 $$
 g\circ h \circ f:
 \begin{cases}
 x\mapsto \phi(\gamma(\alpha(x,u),z),y)\\
 y\mapsto \psi(\gamma(\alpha(x,u),z),y)\\
 z\mapsto \delta(\alpha(x,u),z)\\
u\mapsto \beta(x,u)\\
v\mapsto v
\end{cases} 
 $$ 
This also explains a structure of multiple products.

\sm 
 
  {\bf\punct Proof of Propositions \ref{pr:conjugacy} and \ref{pr:mnogo}.}
The same arguments prove Proposition \ref{pr:mnogo}.
Next, there is a one-to-one correspondence between the sets
$$
\H \setminus\bigl(\G\times \H\bigr)/\H
\,\simeq\,
\G//\H.
$$
The first set is subsemigroup in
$$
\H \setminus\bigl(\G\times \G\bigr)/\H                       
$$                       
and therefore the second set also has a semigroup structure. It remains to verify, that 
two multiplications in   $\G//\H$ coincide.

\sm

{\bf\punct Proof of Theorem \ref{th}.}

\begin{lemma}
{\rm a)} The subspace of $S_\infty[m]$-invariant vectors in $L^2(K^\infty//U)$ coincides with
$H$.

\sm

{\rm b)} The sequence of operators $T(\theta_N)$  converges to the projector $P$
in weak operator topology%
\footnote{See, e.g., \cite{RS}, below in (\ref{eq:long}) we use separate continuity
of product with respect to weak topology.}. 
\end{lemma}

{\sc Proof.} a) For the space $L^2(K^\infty)$ this follows from Hewitt-Savage zero-one low, see, e.g,
\cite{Shi}, Section IV.1, Theorem 3.
The space $L^2(K^\infty//U)$ can be regarded as the space of $U$-invariant functions in 
$L^2(K^\infty)$, 
 and the action of $U$ commutes with action of $S_\infty[m]$.

\sm

b) The statement can be easily proved in a straightforward way.
 However, this is a general fact for continuous
representations of infinite symmetric group (A.Lieberman--G.I.Olshanski, 
see \cite{Olsh-kiado}, \cite{Ner-book}, Section VIII.1, Corollary 5).
\hfill$\square$

\sm

{\sc Proof of Theorem \ref{th}.} Decompose $L^2(K^\infty//U)$ as a direct sum of  $H$
and its orthocomplement.
  Consider an operator
$$
\begin{pmatrix}\ov T(g\circ h)&0\\0&0\end{pmatrix}:L^2(K^\infty//U)\to L^2(K^\infty//U).
$$
It equals $P\,T(g\theta_N h)\, P$ for sufficiently large $N$. Therefore
\begin{multline}
\begin{pmatrix}\ov T(g\circ h)&0\\0&0\end{pmatrix}=
P\, T(g\theta_N h)\, P= P\,T(g\theta_{N+k} h)\, P =\\=
\lim_{j\to\infty}P\, T( g\theta_{j} h)\, P=
\lim_{j\to\infty}P\, T( g)T(\theta_{j}) T( h)\, P
=
\\=
P\, T(g) \Bigl( \lim_{j\to\infty}T(\theta_{j})\Bigr) T(h)\, P= P\,T(g)P T(h)\, P=
\\=
(P\,T(g)\,P) (P\,T( h)\, P)=
\begin{pmatrix}\ov T(g)&0\\0&0\end{pmatrix}\begin{pmatrix}\ov T(h)&0\\0&0\end{pmatrix}
.
\label{eq:long}
\end{multline}
Here $\lim_{j\to\infty}$ is the weak operator limit. Thus we get
$$
\ov T(g\circ h)=\ov T(g) \ov T(h).
$$

{\tt Math.Dept., University of Vienna;

Institute for Theoretical and Experimental Physics (Moscow)

Mech.Math.Dept., Moscow State University,

e-mail: neretin(at) mccme.ru

URL:www.mat.univie.ac.at/$\sim$neretin
}             
                 
\end{document}